\documentclass[12pt,reqno]{amsart}
\usepackage{enumerate, latexsym, amsmath, amsfonts, amssymb, amsthm, color}
\def\pmod #1{\ ({\rm{mod}}\ #1)}

\def\l{\left}
\def\r{\right}
\def\bg{\bigg}
\def\({\bg(}
\def\){\bg)}
\def\t{\text}
\def\f{\frac}

\def\ls{\leqslant}

\def\bi{\binom}

\def\eq{\equiv}

\def\Proof{\noindent{\it Proof}}

\theoremstyle{plain}
\newtheorem{theorem}{Theorem}

\newtheorem{lemma}{Lemma}

\theoremstyle{definition}

\theoremstyle{remark}
\newtheorem{remark}{Remark}

 \vspace{4mm}

\begin{document}
 \baselineskip=17pt
\hbox{Finite Fields Appl. 35(2015), 86--91.}
\medskip

\title
[Proof of a conjectural supercongruence]{Proof of a conjectural supercongruence}

\author
[Xiang-Zi Meng and Zhi-Wei Sun] {Xiang-Zi Meng and Zhi-Wei Sun*}

\thanks{*Supported by the National Natural Science Foundation (grant 11171140)
 of China}

\address {Department of Mathematics, Nanjing
University, Nanjing 210093, People's Republic of China}
\email{{\tt xzmeng@smail.nju.edu.cn}  and  {\tt zwsun@nju.edu.cn}}


 \begin{abstract} Let $m>2$ and $q>0$ be integers with $m$ even or $q$ odd. We show the supercongruence
 $$\sum_{k=0}^{p-1}(-1)^{km}\binom{p/m-q}{k}^m\equiv0\pmod{p^3}$$
for any prime $p>mq$. This confirms a conjecture of Sun.
\end{abstract}

\maketitle

\section{Introduction}
\setcounter{lemma}{0}
\setcounter{theorem}{0}
\setcounter{corollary}{0}
\setcounter{remark}{0}
\setcounter{equation}{0}
\setcounter{conjecture}{0}

Let $p$ be a prime. A $p$-adic congruence is called a {\it supercongruence} if it happens to hold modulo higher powers of $p$. A classical result due to J. Wolstenholme (cf. \cite{W} or \cite{HT}) states that if $p>3$ then
\[\sum_{k=1}^{p-1}\f1k\equiv0\pmod{p^2}\quad \mbox{and} \quad \binom{2p-1}{p-1}\equiv1\pmod{p^3}. \]
For some recent supercongruences modulo $p^2$ involving products of two binomial coefficients, one may consult \cite{S13} and  \cite{GZ}.

Recently, Z.-W. Sun \cite{S15} established some new supercongruences modulo prime powers motivated by the well-known formula
\[\lim_{n\rightarrow\infty}\l(1+\frac1n\r)^n=e.\]
For example, he obtained the following theorem.
\begin{theorem}  For any prime $p>3$,  we have
$$\sum_{k=0}^{p-1}\bi{-1/(p+1)}k^{p+1}\eq0\pmod{p^5}$$
and $$\sum_{k=0}^{p-1}\bi{1/(p-1)}k^{p-1}\eq\f23p^4B_{p-3}\pmod{p^5},$$
where $B_0,B_1,B_2,\ldots$ are the well-known Bernoulli numbers.
\end{theorem}
In this paper we show the following result conjectured by Sun \cite{S15}.
\medskip
\begin{theorem} \label{Th1.2} Let $m>2$ and $q>0$ be integers with $m$ even or $q$ odd. Then, for any prime $p>mq$, we have the supercongruence
\begin{equation}\label{1.1}\sum_{k=0}^{p-1}(-1)^{km}\binom{p/m-q}{k}^m\equiv0\ \pmod{p^3}.
\end{equation}
\end{theorem}
We are going to provide few lemmas in the next section, and then show Theorem 1.2 in Section 3.
\medskip

\maketitle

\section{Some lemmas}
\setcounter{lemma}{0}
\setcounter{theorem}{0}
\setcounter{corollary}{0}
\setcounter{remark}{0}
\setcounter{equation}{0}
\setcounter{conjecture}{0}

\begin{lemma}\label{Lem2.1} For any positive integer $n$, we have
$$\sum_{k=0}^n\bi nk(-1)^k k^m=0\quad\t{for all}\ m=0,1,\ldots,n-1.$$
\end{lemma}
\begin{remark}\label{Rem2.1} This is a well-known result, see, e.g., \cite[pp.\,125-126]{vLW}. \end{remark}

\begin{lemma}\label{Lem2.2} Let $q$ be a positive integer, and let $p>2q$ be a prime. Then, for any integer $k$ with $q-1\ls k\ls p-1$ we have
\begin{equation}\label{2.1}\sum_{q\ls j\ls k}\frac1{j^2}+\sum_{q\ls j<p+q-1-k}\frac1{j^2}\equiv\sum_{j=q}^{p-q}\frac1{j^2}+\sum_{0\ls l<q-1}\frac1{(k-l)^2}\pmod{p}.
\end{equation}
\end{lemma}
\Proof. Observe that
\begin{align*}\sum_{q\ls j<p+q-1-k}\frac1{j^2}\equiv&\sum_{q\ls j<p+q-1-k}\frac1{(p-j)^2}=\sum_{k+1-q<i\ls p-q}\f1{i^2}
\\=&\sum_{0<i<q}\f1{i^2}+\sum_{i=q}^{p-q}\f1{i^2}-\sum_{0<i\ls k+1-q}\f1{i^2}\pmod{p}
\end{align*}
and
$$\sum_{q\ls j\ls k}\frac1{j^2}+\sum_{0<i<q}\f1{i^2}-\sum_{0<i\ls k+1-q}\f1{i^2}=\sum_{k+1-q<j\ls k}\f1{j^2}=\sum_{0\ls l<q-1}\f1{(k-l)^2}.$$
So the desired (\ref{2.1}) follows. \qed

\begin{lemma}\label{Lem2.3} Let $m>2$ and $q>0$ be integers. Then, for any prime $p>mq$ we have
\begin{equation}\label{2.3}\sum_{k=q-1}^{p-1}\binom{k}{q-1}^m\equiv0\pmod{p}\end{equation}
and
\begin{equation}\label{2.4}\sum_{k=q-1}^{p-1}\binom{k}{q-1}^m\sum_{0\ls l<q-1}\frac1{(k-l)^2}\equiv0\pmod{p}.\end{equation}
\end{lemma}
\Proof. As $m(q-1)\ls mq-1\ls p-2$, we can write
$$\bi x{q-1}^m=\sum_{j=0}^{p-2}a_jx^j$$
with $a_0,\ldots,a_{p-2}$ $p$-adic integers. It is well-known that $\sum_{k=0}^{p-1}k^j\eq0\pmod p$ for any positive integer $j\not\eq0\pmod{p-1}$(see, e.g., \cite[p.\, 235]{IR}).
Thus
\begin{align*}\sum_{k=q-1}^{p-1}\bi k{q-1}^m=\sum_{k=0}^{p-1}\sum_{j=0}^{p-2}a_jk^j
=\sum_{j=0}^{p-2}a_j\sum_{k=0}^{p-1}k^j\eq&0\pmod{p}.
\end{align*}

Let $l$ be any integer with $0\ls l<q-1$. Then
\begin{align*} &\sum_{k=q-1}^{p-1}\binom{k}{q-1}^m\frac1{(k-l)^2}
\\\equiv&\sum_{k=l+1}^{p+l-1}\binom{k}{q-1}^m\frac1{(k-l)^2}=\sum_{k=1}^{p-1}\binom{k+l}{q-1}^m\frac1{k^2}\pmod p.
\end{align*}
As $m(q-1)<mq\ls p-1$ and $\bi l{q-1}=0$, we may write
$$\bi{x+l}{q-1}^m=\sum_{m\ls j<p-1}c_jx^j$$
with $c_m,\ldots,c_{p-2}$ $p$-adic integers.
Hence
\begin{align*}\sum_{k=q-1}^{p-1}\binom{k}{q-1}^m\frac1{(k-l)^2}\eq&\sum_{k=0}^{p-1}\sum_{m\ls j<p-1}c_jk^{j-2}
\\\equiv&\sum_{m\ls j<p-1}c_j\sum_{k=0}^{p-1}k^{j-2}\eq0\pmod{p}.
\end{align*}

Now it follows from the above that
\begin{align*}&\sum_{k=q-1}^{p-1}\binom{k}{q-1}^m\sum_{0\ls l<q-1}\frac1{(k-l)^2}
\\=&\sum_{0\ls l<q-1}\sum_{k=q-1}^{p-1}\binom{k}{q-1}^m\frac1{(k-l)^2}\equiv0\pmod{p}.\end{align*}
This concludes the proof. \qed

\section{Proof of Theorem 1.2}

\setcounter{lemma}{0}
\setcounter{theorem}{0}
\setcounter{corollary}{0}
\setcounter{remark}{0}
\setcounter{equation}{0}
\setcounter{conjecture}{0}

\medskip
\noindent{\it Proof of Theorem 1.2}. For each $k\in \{1,2,\cdots,p-1\}$, obviously
\[(-1)^{km}\binom{p/m-q}{k}^m=(-1)^{km}\prod_{j=1}^k\l(\f{p/m-q-j+1}j\r)^m=\prod_{j=1}^k\l(1+\frac{q-1}j-\frac p{jm}\r)^m\]
is congruent to
\begin{align*} &\prod_{j=1}^k\l(\l(1+\frac{q-1}j\r)^m-\l(1+\frac{q-1}j\r)^{m-1}\frac{pm}{jm}+\frac{m(m-1)}2\l(1+\frac{q-1}j\r)^{m-2}\frac{p^2}{j^2m^2}\r)
\\=&\prod_{j=1}^k\l(\l(1+\frac{q-1}j\r)^m-\l(1+\frac{q-1}j\r)^{m-1}\frac{p}{j}+\frac{m-1}{2m}\l(1+\frac{q-1}j\r)^{m-2}\frac{p^2}{j^2}\r)
\end{align*}
module $p^3$, and hence
\begin{align*} &(-1)^{km}\binom{p/m-q}{k}^m-\prod_{j=1}^k\l(\l(1+\frac{q-1}j\r)^m-\l(1+\frac{q-1}j\r)^{m-1}\frac{p}{j}\r)
\\\equiv&\frac{m-1}{2m}\sum_{j=1}^k\l(1+\frac{q-1}j\r)^{m-2}\frac{p^2}{j^2}\prod^k_{i=1\atop i\not=j}\l(1+\frac{q-1}i\r)^m
\\=&p^2\frac{m-1}{2m}\sum_{j=1}^k\frac1{(j+q-1)^2}\prod_{i=1}^k\l(1+\frac{q-1}i\r)^m
\\=&p^2\frac{m-1}{2m}\binom{k+q-1}{k}^m\sum_{j=1}^k\frac1{(j+q-1)^2}\pmod{p^3}.
\end{align*}
Note that
\begin{align*}
&\prod_{j=1}^k\l(\l(1+\frac{q-1}j\r)^m-\l(1+\frac{q-1}j\r)^{m-1}\frac{p}{j}\r)
\\=&\prod_{j=1}^k\l(1+\frac{q-1}j\r)^{m-1}\prod_{j=1}^k\l(1+\frac{q-1}j-\frac pj\r)
\\=&\binom{k+q-1}{k}^{m-1}(-1)^k\binom{p-q}{k}.
\end{align*}
So, from the above we obtain
\begin{align*} &\sum_{k=0}^{p-1}(-1)^{km}\binom{p/m-q}{k}^m-\sum_{k=0}^{p-1}\binom{k+q-1}{k}^{m-1}(-1)^k\binom{p-q}{k}
\\\equiv&\sum_{k=0}^{p-1}p^2\frac{m-1}{2m}\binom{k+q-1}{q-1}^m\sum_{0<j\ls k}\frac1{(j+q-1)^2}
\\=&p^2\frac{m-1}{2m}\sum_{k=q-1}^{p+q-2}\binom{k}{q-1}^m\sum_{q\ls j\ls k}\frac1{j^2}
\\\equiv& p^2\frac{m-1}{2m}\sum_{k=q-1}^{p-1}\binom{k}{q-1}^m\sum_{q\ls j\ls k}\frac1{j^2}\pmod{p^3}.
\end{align*}
Clearly,
\begin{align*}&2\sum_{k=q-1}^{p-1}\bi k{q-1}^m\sum_{q\ls j\ls k}\f1{j^2}
\\=&\sum_{k=q-1}^{p-1}\binom{k}{q-1}^m\sum_{q\ls j\ls k}\frac1{j^2}+\sum_{k=q-1}^{p-1}\binom{p+q-2-k}{q-1}^m\sum_{q\ls j<p+q-1-k}\frac1{j^2},
\end{align*}
and
\[\binom{p+q-2-k}{q-1}^m\equiv\binom{q-2-k}{q-1}^m=(-1)^{m(q-1)}\binom{k}{q-1}^m=\bi k{q-1}^m\pmod{p}\]
for all $k=q-1,q,\ldots,p-1$. Therefore, with the helps of Lemmas 2.2 and 2.3, we get
\begin{align*}&2\sum_{k=q-1}^{p-1}\bi k{q-1}^m\sum_{q\ls j\ls k}\f1{j^2}
\\\eq&\sum_{k=q-1}^{p-1}\binom{k}{q-1}^m\bigg(\sum_{q\ls j\ls k}\frac1{j^2}+\sum_{q\ls j<p+q-1-k}\frac1{j^2}\bigg)
\\\eq&\sum_{k=q-1}^{p-1}\binom{k}{q-1}^m\l(\sum_{j=q}^{p-q}\f1{j^2}+\sum_{0\ls l<q-1}\f1{(k-l)^2}\r)
\\\eq&0\pmod{p}.
\end{align*}
As the degree of the polynomial $\bi{x+q-1}{q-1}^{m-1}$ is $(m-1)(q-1)<(m-1)q<p-q$, by Lemma 2.1 we have
$$\sum_{k=0}^{p-1}\bi{k+q-1}k^{m-1}(-1)^k\bi{p-q}k=\sum_{k=0}^{p-q}\bi{p-q}k(-1)^k\bi{k+q-1}{q-1}^m=0.$$
Therefore the desired (1.1) follows from the above. \qed


\medskip

\begin{thebibliography}{vLW}
\bibitem{GZ} V.J.W. Guo and J. Zeng, {\it Some $q$-analogues of supercongruences of Rodriguez-Villegas}, J. Number Theory {\bf 145} (2014), 301--316.
\bibitem {HT} C. Helou, G. Terjanian, {\it On Wolstenholme's theorem and its converse}, J. Number Theory {\bf 128} (2008), 475--499.
\bibitem {IR} K. Ireland, M. Rosen, A Classical Introduction to Modern Number Theory, 2nd Edition, Grad. Texts in Math., Vol. 84, Springer, New York, 1990.
\bibitem{S13} Z.-W. Sun, {\it Supercongruences involving products of two binomial coefficients}, Finite Fields Appl. {\bf 22} (2013), 24--44.
\bibitem {S15} Z.-W. Sun, {\it Supercongruences motivated by $e$}, J. Number Theory {\bf 147} (2015), 326--341.
\bibitem {vLW} J.H. van Lint, R.M. Wilson, A Course in Combinatorics, 2nd Edition, Cambridge Univ. Press, Cambridge, 2011.
\bibitem {W} J. Wolstenholme, {\it On certain propertities of prime numbers}, Quart. J. Appl. Math. {\bf 5} (1862), 35--39.
\end{thebibliography}
\end{document}